\documentclass[12pt,centertags,oneside]{article}
\usepackage{amsmath,amstext,amsthm,amscd,typearea}
\usepackage{amssymb}
\usepackage{a4wide}
\usepackage[mathscr]{eucal}
\usepackage{mathrsfs}
\usepackage{typearea}
\usepackage{charter}
\usepackage{pdfsync}
\usepackage{url}
\usepackage{bigints}
\usepackage{mathtools}
\usepackage{xcolor}
\usepackage{ulem}

\usepackage[a4paper,width=16.2cm,top=3cm,bottom=3cm]{geometry}

\usepackage{bbm, dsfont}
\usepackage[numbers,sort&compress]{natbib}

\numberwithin{equation}{section}
\usepackage{comment}

\newtheorem{theorem}{Theorem}[section]
\newtheorem{definition}[theorem]{Definition}
\newtheorem{proposition}[theorem]{Proposition}

\newtheorem{lemma}[theorem]{Lemma}
\newtheorem{remark}[theorem]{Remark}

\newtheorem{example}[theorem]{Example}

\newcommand{\supp}{{\rm Supp}}

\newcommand{\vol}{\mathop{\mathrm{vol}}}

\newcommand{\ddc}{dd^c}
\newcommand{\dc}{d^c}

\newcommand{\PSH}{{\rm PSH}}

\newcommand{\capacity}{\mathop{\mathrm{cap}}\nolimits}

\newcommand{\C}{\mathbb{C}}

\newcommand{\N}{\mathbb{N}}

\newcommand{\R}{\mathbb{R}}

\title{\bf Relative non-pluripolar product of currents\\on compact Hermitian manifolds}
\providecommand{\keywords}[1]{\textbf{\textit{Keywords:}} #1}
\providecommand{\subject}[1]{\textbf{\textit{Mathematics Subject Classification 2020:}} #1}

\author{Zhenghao Li and Shuang Su}

\date{\today}
\begin{document}
\maketitle
\begin{abstract}
    Let $(X,\omega)$ be a compact Hermitian manifold. Assume that the Hermitian form $\omega$ satisfies $\partial \overline{\partial} \omega=0$, $\partial \omega \wedge \overline{\partial}\omega=0$. We prove that the relative non-pluripolar product is well-defined on $X$ and satisfies the monotonicity property, generalizing Vu's results in the  K\"ahler setting.
\end{abstract}


\noindent
, {closed positive current}, {Monge--Amp\`{e}re operator}, {relative non-pluripolar product}, {compact Hermitian manifold}. 
\\
\noindent
\subject {32U15}, {32U40}, {53C55}.

\section{Introduction}\label{Section:Introduction}

Pluripotential theory is the study of plurisubharmonic (psh for short) functions and closed positive currents.
One of the central problems in pluripotential theory is to construct a suitable intersection theory of closed positive currents on complex manifolds. More precisely, let $T$ and $S$ be two closed positive currents on a complex manifold $X$, the problem is how to define their intersection $T\wedge S$ in a reasonable way?

When $T$ is of bi-degree $(1,1)$, this problem has been studied extensively. Inspired by the work of Chern--Levine--Nirenberg \cite{Chern-Levine-Nirenberg}, Bedford--Taylor \cite{Bedford-Taylor_Dirichlet, Bedford-Taylor_NewCapacity} gave the first answer to this problem in the local setting. Since $T$ can be written locally as $\ddc u$, where $u$ is a psh function, the problem can be reduced to the study of psh functions. Bedford and Taylor proved that $T\wedge S := \ddc (uS)$ is a closed positive current when $u$ is locally bounded. Later, this construction was extended to certain unbounded potentials with mild singularities (see \cite{Demailly_MongeAmpere, Fornaess-Sibony_OkasInequality}). 

Let $T_{1}, \dots, T_{m}$ be closed positive $(1,1)$-currents on an open subset $U$ of $\C^{n}$. The non-pluripolar product $\langle T_{1} \wedge \dots \wedge T_{m} \rangle$, defined on $U$, was introduced in \cite{Bedford-Taylor_FineTopology}, providing a way to define the wedge product of closed positive $(1,1)$-currents without further assumptions on $T_j$. This construction was later developed on compact K\"ahler manifolds in \cite{Guedj-Zeriahi_WeightedMAEnergy,Boucksom-Eyssidieux-Guedj-Zeriahi}, and has many applications in complex geometry (see \cite{Boucksom-Eyssidieux-Guedj-Zeriahi,Darvas-Xia, Su-Vu_volumepaper,WittNystrom2019Duality}). Let $T$ be a closed positive $(p,p)$-current. In \cite{Vu_RelativeNonPluripolar}, Vu introduced the non-pluripolar product of $T_{1}, \dots, T_{m}$ relative to $T$ on compact K\"ahler manifolds, indicated by
$\langle T_{1} \wedge \dots \wedge T_{m} ~ \dot{\wedge} ~ T \rangle$. This construction generalizes the usual non-pluripolar product. Indeed, when $T=[X]$, it coincides with $\langle T_{1} \wedge \dots \wedge T_{m} \rangle$.

The aim of this paper is to study relative non-pluripolar products on compact complex manifolds that are not necessarily K\"ahler. It was proved in \cite{BottChernvolumes} that, when a compact complex manifold satisfies the bounded volume condition, the non-pluripolar product $\langle T_{1} \wedge \dots \wedge T_{m} \rangle$ is well-defined. In the relative setting, it remains unclear whether the bounded volume condition implies the well-definedness of the corresponding relative non-pluripolar product (see Remark \ref{bvcrelativenon-pluripolar} for more discussion).

In this paper, we consider the compact Hermitian manifold $(X,\omega)$ satisfying
\begin{equation}
\label{maincondition}
     \partial \overline{\partial} \omega = 0,\quad \partial \omega \wedge \overline{\partial} \omega = 0.
\end{equation}
The following is our first main result, which generalizes \cite[Lemma 3.4]{Vu_RelativeNonPluripolar} from the K\"ahler setting to the Hermitian setting.
\begin{theorem}[Theorem \ref{Theorem:WellDefined}]
    Let $(X,\omega)$ be an $n$-dimensional compact Hermitian manifold such that $\omega$ satisfies (\ref{maincondition}). Let $T_j$ be a closed positive $(1,1)$-current on $X$ for $j=1,\dots,m$. Let $T$ be a closed positive $(p,p)$-current on $X$ such that $m+p \leq n$. Then the relative non-pluripolar product $\langle T_{1} \wedge \dots \wedge T_{m} ~\dot{\wedge}~ T \rangle$ is well-defined.
\end{theorem}
Although condition \eqref{maincondition} is stronger than the bounded volume condition, it guarantees that the relative non-pluripolar product is well-defined for every closed positive $(p,p)$-current $T$. In Example \ref{Example:nonKahler}, we give explicit examples of non-K\"ahler manifolds of dimension at least $2$ that satisfy (\ref{maincondition}). 

Let $T$ be a closed $(p,p)$-current. We denote by $\{T\}_{BC} \in H^{p,p}_{BC}(X,\C)$ its Bott--Chern class. If $T$ is positive, then it is real and hence $\{T\}_{BC} \in H^{p,p}_{BC}(X,\R)$. A class $\alpha\in H^{p,p}_{BC}(X,\mathbb{R})$ is called \textit{pseudoeffective} if it can be represented by a closed positive $(p,p)$-current.
Let $\theta$ be a real closed $(1,1)$-form on $X$ such that its Bott--Chern class $\{\theta\}_{BC}$ is pseudoeffective. Let $T,T'\in \{\theta\}_{BC}$ be two closed positive $(1,1)$-currents. We can write
\[
    T= \ddc \tilde{u} + \theta, \quad T' = \ddc \tilde{u}' +\theta,
\]
where $\tilde{u},\tilde{u}' \in \PSH(X,\theta)$. We say $T'$ is \textit{less singular (more singular)} than $T$ if $\tilde{u} \leq \tilde{u}' + \mathcal{O}(1)$ ($\tilde{u}' \leq \tilde{u} + \mathcal{O}(1)$). Our second main result is the monotonicity property of the relative non-pluripolar product in the Hermitian setting.

\begin{theorem}[Theorem~\ref{hermitianmonotonicity}]
\label{hermitianmonotonicityintro}
    Let $(X,\omega)$ be an $n$-dimensional compact Hermitian manifold such that $\omega$ satisfies (\ref{maincondition}). For $j=1,\dots, m$, let $\{\theta_{j}\}_{BC}$ be a pseudoeffective Bott--Chern class, and $T_{j},T'_{j}  \in \{\theta_{j}\}_{BC}$ be closed positive $(1,1)$-currents such that $T'_{j}$ is less singular than $T_{j}$. Let $T$ be a closed positive $(p,p)$-current on $X$ such that $m+ p \leq n$. Then we have  
    \[
        \int_{X} \langle T_{1} \wedge \dots\wedge T_{m} ~\dot{\wedge}~ T \rangle \wedge \omega^{n-m-p} \leq \int_{X}\langle T'_{1} \wedge \dots \wedge T'_{m} ~\dot{\wedge}~ T \rangle \wedge \omega^{n-m-p}.
    \]
\end{theorem}

Theorem \ref{hermitianmonotonicityintro} extends Vu's monotonicity theorem for relative non-pluripolar products on compact K\"ahler manifolds \cite[Theorem~1.1]{Vu_RelativeNonPluripolar} to the Hermitian setting. We refer the reader to \cite{Boucksom-Eyssidieux-Guedj-Zeriahi,Darvas-DiNezza-Lu_Monotonicity,WittNystrom_MonotonicityNonPluripolarMAMasses} for related monotonicity results for non-pluripolar products on compact K\"ahler manifolds, and to \cite{BottChernvolumes} for corresponding results on compact complex manifolds. We also note that, under condition \eqref{maincondition}, the monotonicity property above holds in the sense of cohomology classes, as in \cite{Vu_RelativeNonPluripolar} (see Remark~\ref{monotonicityclassesremark}).

The organization of this paper is as follows. In Section 2, we recall the construction of relative non-pluripolar products and review some of their properties. In Section 3, we prove our two main results.

\subsection*{Acknowledgement}
The authors thank Duc-Viet Vu for posting this problem. The authors also thank Yinji Li, Zhuo Liu, and Ngoc Cuong Nguyen for discussions and valuable comments. Zhenghao Li would like to express his sincere gratitude to Xiangyu Zhou for his encouragement and support. He also thanks the members of the Complex Analysis Group at the University of Cologne for their warm hospitality during his visit.

The research of Zhenghao Li and Shuang Su was partially supported by the Institute for Basic Science (IBS-R032-D1-2026-a00). Shuang Su was also partially supported by the Deutsche Forschungsgemeinschaft (Project number 500055552) and the ANR--DFG grant QuaSiDy (ANR-21-CE40-0016) during his Ph.D. studies.

\section{Construction of relative non-pluripolar products}
\label{Section:RelativeNonPluripolar}
In this section, we recall the construction of relative non-pluripolar products and some of their properties, following the approach of \cite{Vu_RelativeNonPluripolar} while providing more details in the proofs for the reader's convenience. 

We first recall some notions concerning currents that will be used later. Let $(X,\omega)$ be an $n$-dimensional Hermitian manifold and $T$ be a positive $(p,p)$-current on $X$. The trace measure of $T$ with respect to $\omega$ is defined by $T\wedge\omega^{n-p}$, which is a positive Radon measure on $X$. When the choice of $\omega$ is clear from the context, we denote 
\[
    \lVert T\rVert_K:=\int_K T\wedge \omega^{n-p}
\]
for every Borel subset $K$ of $X$ and denote $\lVert T\rVert:=\lVert T\rVert_X$ for simplicity. If $K\Subset X$, then $\lVert T\rVert_K<\infty$. Throughout this section, we will frequently consider a sequence of positive currents whose trace measures are uniformly bounded on every relatively compact subset of $X$. This property is independent of the choice of the Hermitian form $\omega$.

Let $U\subset\C^n$ be an open subset, $\omega:=i\Sigma_{j=1}^n dz_j\wedge d\bar{z}_j$ be the standard K\"ahler form on $\C^n$, and $K\Subset U$ be a Borel  subset. Let $T$ be a $(p,p)$-current of order $0$ on $U$. The mass of $T$ on $K$ is defined by
\begin{align*}
    \lvert T \rvert_K:=\sup_{\Phi}\lvert \langle T,\Phi\rangle \rvert,
\end{align*}
where the supremum is taken over all $(n-p,n-p)$-forms $\Phi$ on $U$ such that $\supp\Phi\subset K$ and the absolute value of their coefficients are bounded by $1$. If $T$ is positive, then there exists a constant $c_{n,p}>0$ that depends only on $n$ and $p$, such that
\begin{align*}
    c_{n,p}^{-1}~\lvert T \rvert_K\leq \lVert T\rVert_K=\int_K T\wedge\omega^{n-p} \leq c_{n,p}~\lvert T\rvert_K.    
\end{align*}

We next introduce a property of increasing sequences of positive currents, which will be used in the construction of relative non-pluripolar products. 
\begin{proposition}
\label{Proposition:IncreasingBoundedMeasures}
    Let $U \subset \C^{n}$ be an open subset. Assume that $(T_k)_{k\in\N}$ is an increasing sequence of positive $(p,p)$-currents on $U$, and
    \[
        \sup_{k\in\N}\|T_k\|_K <\infty
    \]
    for every compact subset $K$ of $U$. Then $T_k$ weakly converges to a positive current $T$, and
    \[
        \langle T,\Phi \rangle = \lim_{k\to\infty} \langle T_k,\Phi\rangle
    \]
    for every bounded Borel $(n-p,n-p)$-form $\Phi$ on $U$ with compact support.
\end{proposition}

\begin{proof}
    For every compact subset $K$ of $U$, the sequence $(\lVert T_k\rVert_K)_{k\in\N}$ is increasing and bounded, hence its limit exists. Write
    \begin{align*}
        T_k=i^{p^2}\sum_{|J|=|L|=p}T_{k,J,L}\ dz_J\wedge d \bar{z}_L,
    \end{align*}
    where $T_{k,J,L}$ are complex Radon measures. For any multi-indices $J,L$ and $\varphi\in \mathscr{C}_c^\infty(U,\C)$, 
	\begin{align*}
	    \lvert\langle T_{k,J,L},\varphi\rangle\rvert\leq \sup_U |\varphi|\cdot\lvert T_{k}\rvert_{\supp\varphi}\leq c_{n,p}\sup_U |\varphi|\cdot\lVert T_k \rVert_{\supp\varphi}.
	\end{align*}
    Hence, the sequence $(\langle T_{k,J,L},\varphi\rangle)_{k\in\N}$ is bounded. Then, to prove $\lim_{k\to\infty}\langle T_{k,J,L},\varphi\rangle$ exists, it suffices to show that any convergent subsequences $(\langle T_{j_k,J,L},\varphi\rangle)_{k\in\N}$ and $(\langle T_{\ell_k,J,L},\varphi\rangle)_{k\in\N}$ have the same limit. To prove this, we can assume that $j_k< \ell_k$ for every $k\in\mathbb{N}$. Then
    \begin{align*}
        \lvert\langle T_{\ell_k,J,L} - T_{j_k,J,L},\varphi\rangle\rvert\leq \sup_U |\varphi|\cdot\lvert T_{\ell_k}-T_{j_k}\rvert_{\supp\varphi}\leq c_{n,p}\sup_U |\varphi|\cdot\lVert T_{\ell_k}-T_{j_k} \rVert_{\supp\varphi},
    \end{align*}
	which tends to $0$ as $k\to\infty$ because the sequence $(\lVert T_k\rVert_{\supp\varphi})_{k\in\N}$ is convergent. It follows that $(\langle T_{j_k,J,L},\varphi\rangle)_{k\in\N}$ and $(\langle T_{\ell_k,J,L},\varphi\rangle)_{k\in\N}$ have the same limit, and hence $\lim_{k\to\infty}\langle T_{k,J,L},\varphi\rangle$ exists. Therefore, $(T_k)_{k\in\N}$ weakly converges to a positive $(p,p)$-current $T$ on $U$. Moreover, for every strongly positive smooth $(n-p,n-p)$-form $\Phi$ on $U$ with compact support, $(\langle T_k,\Phi\rangle)_{k\in\N}$ is increasing and converges to $\langle T,\Phi\rangle$, hence $T_k\leq T$ for every $k\in\N$.
	
	Let $\Phi$ be a bounded Borel $(n-p,n-p)$-form on $U$ with compact support such that the absolute values of its coefficients are bounded by a constant $c_\Phi>0$. By Lusin's Theorem (see, e.g., \cite[2.24]{Rudin_RealComplexAnalysis}), for every $\varepsilon>0$ there exists a continuous $(n-p,n-p)$-form $\Phi_{\varepsilon}$ on $U$ with compact support such that
	\[
	      \|T\|_{\{x\in U:\Phi(x)\neq\Phi_{\varepsilon}(x)\}}\leq\varepsilon
	\]
    and the absolute values of the coefficients of $\Phi_{\varepsilon}$ are bounded by $c_\Phi$.
	Since $\langle T_k,\Phi_{\varepsilon}\rangle$ converges to $\langle T,\Phi_{\varepsilon}\rangle$, there exists $N\in\mathbb{N}$ such that 
	\[
	    |\langle T_k,\Phi_{\varepsilon}\rangle-\langle T,\Phi_{\varepsilon}\rangle|\leq\varepsilon.
	\]
    for every $k\geq N$.
	It follows that for every $k\geq\N$ we have
    \begin{align*}
        |\langle T_k,\Phi \rangle - \langle T,\Phi \rangle |
        &\leq |\langle T_k,\Phi \rangle - \langle T_k,\Phi_{\varepsilon} \rangle | + |\langle T_k,\Phi_{\varepsilon} \rangle - \langle T,\Phi_{\varepsilon} \rangle |+ |\langle T,\Phi_{\varepsilon} \rangle - \langle T,\Phi \rangle |\\
        &\leq 2 c_\Phi |T_k|_{\{x\in U:\Phi(x)\neq\Phi_{\varepsilon}(x)\}} +\varepsilon+2 c_\Phi |T|_{\{x\in U:\Phi(x)\neq\Phi_{\varepsilon}(x)\}}\\
        &\leq 2 c_\Phi c_{n,p}\lVert T_k\rVert_{\{x\in U:\Phi(x)\neq\Phi_{\varepsilon}(x)\}}+\varepsilon+ 2 c_\Phi c_{n,p}\lVert T\rVert_{\{x\in U:\Phi(x)\neq\Phi_{\varepsilon}(x)\}}\\
        &\leq 4 c_\Phi c_{n,p} \lVert T\rVert_{\{x\in U:\Phi(x)\neq\Phi_{\varepsilon}(x)\}} +\varepsilon\\
        &\leq 4c_\Phi c_{n,p}\varepsilon +\varepsilon.
    \end{align*}
    Therefore, we have $\langle T,\Phi \rangle = \lim_{k\to\infty} \langle T_k,\Phi\rangle$.
\end{proof}

The following property of Monge--Amp\`ere operators was first established in \cite[Corollary~4.3]{Bedford-Taylor_FineTopology} and later generalized to the relative setting in \cite{Vu_RelativeNonPluripolar}.

\begin{theorem}(\cite[Theorem 2.9]{Vu_RelativeNonPluripolar})
\label{Theorem:PlurifinelyTopologicalProperty}
    Let $U\subset\C^n$ be an open subset. Let $T$ be a closed positive current on $U$. Let $u_j, u'_j$ be locally bounded psh functions on $U$ for $j=1, \dots, m$, and $v_k, v'_k$ be locally bounded psh functions on $U$ for $k=1,\dots,\ell$. Assume that $u_j=u'_j$ on $W:=\cap_{k=1}^{\ell}\{v_k>v'_k\}$ for $j=1,\dots,m$. Then we have
    \[
        \mathbbm{1}_W \ddc u_1 \wedge \dots \wedge \ddc u_m \wedge T=\mathbbm{1}_W \ddc u'_1 \wedge \dots \wedge \ddc u'_m \wedge T.
    \]
\end{theorem}

\begin{remark}
\label{Remark: Plurifine topological property}
    Recall that a quasi-plurisubharmonic (quasi-psh for short) function can be locally expressed as the sum of a psh function and a smooth function. For any finite collection of quasi-psh functions, one can locally find a common smooth function such that adding it to each quasi-psh function yields a psh function. Therefore, as pointed out in \cite[Remark 2.10]{Vu_RelativeNonPluripolar}, Theorem \ref{Theorem:PlurifinelyTopologicalProperty} remains valid if the locally bounded psh functions in the theorem are replaced by locally bounded quasi-psh functions and $U$ is replaced by any complex manifold.
\end{remark}

Now, let $X$ be an $n$-dimensional complex manifold,  $T$ be a closed positive $(p,p)$-current on $X$, and $T_1,\dots,T_m$ be closed positive $(1,1)$-currents on X such that $m+p\leq n$.

We first construct relative non-pluripolar product in the local setting. Let $U$ be a local chart of $X$ and $u_1,\dots,u_m\in PSH (U)$ such that $T_j=dd^c u_j$ on $U$ for $j=1,\dots,m$. For every $k\in\N$, set 
\[
    u_{j,k}:=\max\{u_j,-k\},\quad R_k:=\ddc u_{1,k}\wedge\dots\wedge\ddc u_{m,k}\wedge T|_U
\]
Then $u_{j,k}$ is a locally bounded psh function on $U$, and $R_k$ is a closed positive $(m+p,m+p)$-current on $U$. Since $\{u_j>-k\}=\{u_{j,k}>-k\}$, by Theorem \ref{Theorem:PlurifinelyTopologicalProperty} we have 
\begin{align}
\label{Equation:PlurifinelyTopology_Psh}
    \mathbbm{1}_{\cap_{j=1}^{m} \{ u_j > -k\} }R_k=\mathbbm{1}_{\cap_{j=1}^{m} \{ u_j > -k\} }R_{\ell}\leq\mathbbm{1}_{\cap_{j=1}^{m} \{ u_j > -\ell\} }R_\ell    
\end{align}
whenever $\ell> k$. The sequence of positive currents $(\mathbbm{1}_{\cap_{j=1}^{m} \{ u_j > -k\} }R_k)_{k\in\N}$ is increasing. This immediately yields the following result.

\begin{lemma}(\cite[Lemma 3.1]{Vu_RelativeNonPluripolar})
\label{Lemma:UniformlyBounded_Psh}
    Assume that
    \begin{equation}
    \label{Equation:UniformlyBounded_Psh}
        \sup_{k\in\N} \| \mathbbm{1}_{\cap_{j=1}^{m} \{ u_j > -k\} }R_k\|_K<\infty  
    \end{equation}
    for every compact subset $K$ of $U$. Then the following statements hold.
    \begin{itemize}
        \item[\upshape(i)]
        $\mathbbm{1}_{\cap_{j=1}^{m} \{ u_j > -k\} }R_k$ weakly converges to a positive current $R$ on $U$. Moreover,
        \[
            \langle R,\Phi \rangle = \lim_{k\to\infty} \langle \mathbbm{1}_{\cap_{j=1}^{m} \{ u_j > -k\} }R_k,\Phi\rangle
        \]
        for every bounded Borel $(n-m-p,n-m-p)$-form $\Phi$ on $U$ with compact support.
        
        \item[\upshape(ii)]
        $\mathbbm{1}_{\cap_{j=1}^{m} \{ u_j > -k\} }R=\mathbbm{1}_{\cap_{j=1}^{m} \{ u_j > -k\} }R_k$ for every $k\in\N$, and $\mathbbm{1}_{\cup_{j=1}^m\{u_j =-\infty\}}R=0$. 
    \end{itemize}
\end{lemma}

\begin{proof}
    (i) is deduced by Proposition \ref{Proposition:IncreasingBoundedMeasures}. For any smooth $(n-m-p,n-m-p)$-form $\Phi$ on $U$ with compact support, by (i) we have
    \begin{align*}
         &\langle \mathbbm{1}_{\cap_{j=1}^{m} \{ u_j > -k\} }R, \Phi \rangle=\lim_{\ell\to\infty} \langle \mathbbm{1}_{\cap_{j=1}^{m} \{ u_j > -\ell\} }R_\ell, \mathbbm{1}_{\cap_{j=1}^{m} \{ u_j > -k\} } \Phi \rangle=\langle\mathbbm{1}_{\cap_{j=1}^{m} \{ u_j > -k\} }R_k, \Phi\rangle,\\
         &\langle \mathbbm{1}_{\cup_{j=1}^m\{u_j =-\infty\}}R, \Phi \rangle=\lim_{\ell\to\infty}\langle \mathbbm{1}_{\cap_{j=1}^{m} \{ u_j > -\ell\} }R_\ell, \mathbbm{1}_{\cup_{j=1}^m\{u_j =-\infty\}}\Phi\rangle=0.
    \end{align*}
    This finishes the proof.
\end{proof}

To construct relative non-pluripolar products in the global setting, we introduce the following global analogue of Lemma \ref{Lemma:UniformlyBounded_Psh}, where psh functions are replaced by quasi-psh functions. For $j=1,\dots,m$, let $\theta_j$ be a real closed $(1,1)$-form on $X$ and $\tilde{u}_j\in PSH (X,\theta_j)$ such that $T_j=dd^c \tilde{u}_j + \theta_j$ on $X$. For every $k\in\N$, set
\[
    \tilde{u}_{j,k}:=\max\{\tilde{u}_j,-k\},\quad \tilde{R}_k:=(\ddc \tilde{u}_{1,k}+\theta_1)\wedge\dots\wedge(\ddc \tilde{u}_{m,k}+\theta_m)\wedge T.
\]
Then $\tilde{u}_{j,k}$ is a locally bounded quasi-psh function on $X$, and $\tilde{R}_k$ is a $(m+p,m+p)$-current of order $0$ on $X$ but may not be positive.  Since $\{\tilde{u}_j>-k\}=\{\tilde{u}_{j,k}>-k\}$, by Theorem \ref{Theorem:PlurifinelyTopologicalProperty} and Remark \ref{Remark: Plurifine topological property} we have
\begin{equation}
\label{Equation:PlurifinelyTopology_QuasiPsh}
    \mathbbm{1}_{\cap_{j=1}^{m} \{ \tilde{u}_j > -k\} }\tilde{R}_k=\mathbbm{1}_{\cap_{j=1}^{m} \{ \tilde{u}_j > -k\} }\tilde{R}_\ell
\end{equation}
whenever $\ell> k$. Together with Lemma \ref{Lemma:UniformlyBounded_Psh}, this yields the following result.

\begin{lemma}(\cite[Lemma 3.2]{Vu_RelativeNonPluripolar})
\label{Lemma:UniformlyBounded_QuasiPsh}
    The following statements hold.
    \begin{itemize}
        \item[{\upshape(i)}]
        Foe every $k\in\N$, $\mathbbm{1}_{\cap_{j=1}^{m} \{ \tilde{u}_j > -k\} }\tilde{R}_k$ is a positive current on $X$. 
        
        \item[{\upshape(ii)}]
        Assume that
        \begin{equation}
        \label{Equation:UniformlyBounded_QuasiPsh}
            \sup_{k\in\N} \| \mathbbm{1}_{\cap_{j=1}^{m} \{ \tilde{u}_j > -k\} }\tilde{R}_k\|_K<\infty.
        \end{equation}
        for every compact subset $K$ of $X$. Then $\mathbbm{1}_{\cap_{j=1}^{m} \{\tilde{u}_j > -k\} }\tilde{R}_k$ weakly converges to a positive current $\tilde{R}$ on $X$.

        \item[\upshape(iii)]
        $\mathbbm{1}_{\cap_{j=1}^{m} \{ \tilde{u}_j > -k\} }\tilde{R}=\mathbbm{1}_{\cap_{j=1}^{m} \{ \tilde{u}_j > -k\} }\tilde{R}_k$ for every $k\in\N$, and $\mathbbm{1}_{\cup_{j=1}^m\{\tilde{u}_j =-\infty\}}\tilde{R}=0$.
    \end{itemize}
\end{lemma}

\begin{proof}
    (i) For every $x\in X$, let $U$ be a local chart containing $x$. By the local $dd^c$-lemma, after shrinking $U$ we can find $h_j\in\mathscr{C}^\infty(U,\R)$ for $j=1,\dots,m$ such that $dd^c h_j=\theta_j|_U$. Hence, $u_j:=\tilde{u}_j|_U +h_j\in PSH (U)$ and $T_j=dd^c u_j$ on $U$. We may assume that $h_j$ is bounded on $U$ after further shrinking $U$. Let $c$ be an integer greater than $\sum_{j=1}^m \sup_U |h_j|$. For every $k\in\N$, set
    \[
        u_{j,k}:=\max\{u_j,-k\},\quad R_k:=\ddc u_{1,k}\wedge\dots\wedge\ddc u_{m,k}\wedge T|_U.
    \]
    Since $\{\tilde{u}_j|_U>-k\}\subset \{u_j>-(k+c)\}$ for every $k\in\N$, on $\{\tilde{u}_j|_U>-k\}$ we have
    \[
    \tilde{u}_{j,k+c}=\tilde{u}_j=u_j-h_j=u_{j,k+c}-h_j.
    \]
    Then it follows from (\ref{Equation:PlurifinelyTopology_QuasiPsh}), Theorem \ref{Theorem:PlurifinelyTopologicalProperty}, and Remark \ref{Remark: Plurifine topological property} that
    \begin{align}
    \label{Equation: relationship between Rk and tilteRk}
        \begin{aligned}
            &\left(\mathbbm{1}_{\cap_{j=1}^m\{\tilde{u}_j>-k\}}\tilde{R}_k\right)\big|_U\\
            &=\left(\mathbbm{1}_{\cap_{j=1}^m\{\tilde{u}_j>-k\}}\tilde{R}_{k+c}\right)\big|_U\\
            &=\mathbbm{1}_{\cap_{j=1}^m\{\tilde{u}_j|_U>-k\}} dd^c \left(\tilde{u}_{1,k+c}|_U+\theta_1|_U\right)\wedge\dots\wedge dd^c \left(\tilde{u}_{m,k+c}|_U+\theta_m|_U\right)\wedge T|_U\\
            &=\mathbbm{1}_{\cap_{j=1}^m\{\tilde{u}_j|_U>-k\}} dd^c u_{1,k+c}\wedge\dots\wedge dd^c u_{m,k+c}\wedge T|_U\\
            &=\mathbbm{1}_{\cap_{j=1}^m\{\tilde{u}_j|_U>-k\}} R_{k+c},
        \end{aligned}
    \end{align}
    which is a positive current on $U$. Therefore, $\mathbbm{1}_{\cap_{j=1}^m\{\tilde{u}_j>-k\}}\tilde{R}_k$ is a positive current on $X$. 

    (ii) We retain the notation from (i). For every $k\in\N$, we have
    \[
        \{u_j>-k\}\subset \{\tilde{u}_j|_U>-(k+c)\}\subset\{u_j>-(k+2c)\}.
    \]
    Then it follows from Theorem \ref{Theorem:PlurifinelyTopologicalProperty} that
    \begin{align*}
        \mathbbm{1}_{\cap_{j=1}^m\{u_j>-k\}}R_k=\mathbbm{1}_{\cap_{j=1}^m\{u_j>-k\}}R_{k+2c}\leq\mathbbm{1}_{\cap_{j=1}^m\{\tilde{u}_j|_U>-(k+c)\}} R_{k+2c}\leq\mathbbm{1}_{\cap_{j=1}^m\{u_j>-(k+2c)\}}R_{k+2c}.
    \end{align*}
    Combining this with (\ref{Equation: relationship between Rk and tilteRk}), we deduce that
    \begin{equation}
    \label{Equation:Inequality}
        \mathbbm{1}_{\cap_{j=1}^m\{u_j>-k\}}R_k\leq\left(\mathbbm{1}_{\cap_{j=1}^m\{\tilde{u}_j>-(k+c)\}}\tilde{R}_{k+c}\right)\big|_U\leq\mathbbm{1}_{\cap_{j=1}^m\{u_j>-(k+2c)\}}R_{k+2c}.
    \end{equation}
    It follows from the first inequality in (\ref{Equation:Inequality}) that $(\mathbbm{1}_{\cap_{j=1}^m\{u_j>-k\}}R_k)_{k\in\N}$ satisfies (\ref{Equation:UniformlyBounded_Psh}) and hence by Lemma \ref{Lemma:UniformlyBounded_Psh} (i) it weakly converges to a positive current $R$ on $U$. By (\ref{Equation:Inequality}) again,
    \begin{align*}
        \lim_{k\to\infty} \langle \mathbbm{1}_{\cap_{j=1}^m\{\tilde{u}_j>-k\}}\tilde{R}_k , \Phi \rangle = \langle R, \Phi\rangle
    \end{align*}
    for every strongly positive smooth $(n-m-p,n-m-p)$-form $\Phi$ on $U$ with compact support. Since every test form is a linear combination of strongly positive test forms, the above limit holds for every smooth $(n-m-p,n-m-p)$-form $\Phi$ on $U$ with compact support.
    Therefore, we conclude that $\mathbbm{1}_{\cap_{j=1}^{m} \{\tilde{u}_j > -k\} }\tilde{R}_k$ weakly converges to a positive current $\tilde{R}$ on $X$. Moreover, $\tilde{R}=R$ on $U$.

    (iii) We retain the notation from (i). It follows from Lemma \ref{Lemma:UniformlyBounded_Psh} (ii) and (\ref{Equation: relationship between Rk and tilteRk}) that
    \begin{align*}
        \left(\mathbbm{1}_{\cap_{j=1}^m\{\tilde{u}_j>-k\}}\tilde{R}\right)\big|_U&=\mathbbm{1}_{\cap_{j=1}^m\{\tilde{u}_j|_U>-k\}} R=\mathbbm{1}_{\cap_{j=1}^m\{\tilde{u}_j|_U>-k\}}\mathbbm{1}_{\cap_{j=1}^m\{u_j>-(k+c)\}} R\\
        &=\mathbbm{1}_{\cap_{j=1}^m\{\tilde{u}_j|_U>-k\}}\mathbbm{1}_{\cap_{j=1}^m\{u_j>-(k+c)\}} R_{k+c}=\mathbbm{1}_{\cap_{j=1}^m\{\tilde{u}_j|_U>-k\}}R_{k+c}\\
        &=\left(\mathbbm{1}_{\cap_{j=1}^m\{\tilde{u}_j>-k\}}\tilde{R}_k\right)\big|_U,\\
        \left(\mathbbm{1}_{\cup_{j=1}^m\{\tilde{u}_j =-\infty\}}\tilde{R}\right)\big|_U&=\mathbbm{1}_{\cup_{j=1}^m\{u_j =-\infty\}}R=0.
    \end{align*}
    This finishes the proof.
\end{proof}

\begin{remark}
    The proof of Lemma \ref{Lemma:UniformlyBounded_QuasiPsh} shows that both (\ref{Equation:UniformlyBounded_QuasiPsh}) and the limit current $\tilde{R}$ are independent of the auxiliary choices of $\theta_j$ and $\tilde{u}_j$. As a consequence, in Lemma \ref{Lemma:UniformlyBounded_Psh}, both (\ref{Equation:UniformlyBounded_Psh}) and the limit current $R$ are independent of the auxiliary choices of $u_j$.

    Assume that (\ref{Equation:UniformlyBounded_QuasiPsh}) holds. Then (\ref{Equation:UniformlyBounded_Psh}) holds on every local chart $U$ of $X$, and $\tilde{R}=R$ on $U$. Conversely, assume that there exists a family of local charts covering X such that \eqref{Equation:UniformlyBounded_Psh} holds on each chart. Then \eqref{Equation:UniformlyBounded_QuasiPsh} holds, and $\tilde{R}=R$ on every chart in the cover.
\end{remark}

We are now ready to define the relative non-pluripolar product.
\begin{definition}
    We say that the non-pluripolar product relative to $T$ of $T_1,\dots,T_m$ is well-defined if (\ref{Equation:UniformlyBounded_QuasiPsh}) holds, or equivalently, (\ref{Equation:UniformlyBounded_Psh}) holds for a family of local charts covering X.
    In this case, the non-pluripolar product relative to $T$ of $T_1,\dots,T_m$, which is denoted by $\langle T_{1} \wedge \dots\wedge T_{m} ~\dot{\wedge}~ T \rangle$, is defined as the current $\tilde{R}$ in Lemma \ref{Lemma:UniformlyBounded_QuasiPsh}.
\end{definition}

We refer the reader to \cite[Proposition 3.5 and Proposition 3.6]{Vu_RelativeNonPluripolar} for some properties of relative non-pluripolar products.

The following theorem, concerning the closedness of relative non-pluripolar products, was proved in \cite{Vu_RelativeNonPluripolar}. It generalizes the closedness result for the usual non-pluripolar product \cite[Theorem~1.8]{Boucksom-Eyssidieux-Guedj-Zeriahi}. We provide the proof here for the reader's convenience.

\begin{theorem}(\cite[Theorem 3.7]{Vu_RelativeNonPluripolar})
    $\langle T_1\wedge\dots\wedge T_m ~\dot{\wedge}~ T\rangle$ is closed if it is well-defined.
\end{theorem}

\begin{proof}
    For every $x\in X$, we will show the closedness of the restriction of $\langle T_1\wedge\dots\wedge T_m ~\dot{\wedge}~ T\rangle$ on an open neighborhood of $x$. Let $U$ be a local chart of $X$ containing $x$ and $u_j\in PSH (U)$ for $j=1,\dots,m$, such that $u_j\leq 0$ and $T_j=dd^c u_j$ on $U$. For every $k\in\N$, set 
    \[
        u_{j,k}:=\max\{u_j,-k\},\quad R_k:=\ddc u_{1,k}\wedge\dots\wedge\ddc u_{m,k}\wedge T|_U,\quad R:=\lim_{k\to\infty} \mathbbm{1}_{\cap_{j=1}^{m} \{ u_j > -k\} }R_k.
    \]
    To prove that \(R\) is closed, it suffices to show that $d(\rho \wedge R) = 0$ for every strongly positive \((n - m - p - 1, n - m - p - 1)\)-form $\rho$ on $U$ with constant coefficients.

    For every $k\in\N$, set  
    \[
        \frac{1}{k} \max\left\{ \sum_{j=1}^m u_j, -k \right\} + 1=:\psi_k \in PSH (U).
    \]  
    Then $0\leq\psi_k\leq 1$, $\psi_k=0$ on $\cup_{j=1}^m\{u_j\leq -k\}$, and $\psi_k$ increases to \(\mathbbm{1}_{\cap_{j=1}^m \{u_j > -\infty\}}\). Take $g\in \mathscr{C}^\infty(\R,\R)$ such that $g,g'\geq 0$, $g(0)=0$, $g(1)=1$, $g'(0)=g'(1)=0$, and $g'$ is bounded.
    
    By Lemma \ref{Lemma:UniformlyBounded_Psh} (ii) and the properties of $g$, we have
    \begin{align*}
        &R=\mathbbm{1}_{\cap_{j=1}^m\{u_j>-\infty\}}R=\lim_{k\to\infty}g(\psi_k)R,\\
        &g(\psi_k)R=g(\psi_k)\mathbbm{1}_{\cap_{j=1}^m\{u_j>-k\}}R=g(\psi_k)\mathbbm{1}_{\cap_{j=1}^m\{u_j>-k\}}R_k=g(\psi_k)R_k,
    \end{align*}
    which implies that
    \[
        d(\rho\wedge R)=\lim_{k\to\infty}d\left(g(\psi_k)\rho\wedge R_k\right).
    \]
    Since $\rho\wedge R_k$ is a closed positive current on $U$, it follows from \cite[Lemma 1.9]{Boucksom-Eyssidieux-Guedj-Zeriahi} that
    \[
        d(g(\psi_k)\, \rho \wedge R_k) = g'(\psi_k) d \psi_{k} \wedge \rho \wedge R_k.
    \]
    Therefore, we only need to show that 
    \[
        \lim_{k\to\infty}\int_U g'(\psi_k) \eta\wedge d \psi_{k} \wedge \rho \wedge R_k =0
    \]
    for every smooth $1$-form $\eta$ on $U$ with compact support.

    Write $\eta=\eta_{1}+ \overline{\eta_{2}}$, where $\eta_{1}$ and $\eta_{2}$ are $(1,0)$-forms on $U$.
    Take $\tau\in\mathscr{C}_c^\infty(U,\R)$ such that $\tau=1$ on an open neighborhood of $\supp\eta$. By applying the Cauchy–Schwarz inequality to the standard regularization of $\psi_k$ and then passing to the limit, we have
    \begin{equation}
    \label{closdnesseq0}
        \begin{aligned}
             &\left| \int_{U} g'(\psi_k) \eta \wedge d\psi_k \wedge  \rho \wedge R_k \right|^2 \\
        & \lesssim\left(\int_{U} \tau^2 d\psi_k \wedge d^c \psi_k \wedge \rho \wedge R_k\right) \left(\int_U i(g'(\psi_k))^2 (\eta_{1} \wedge \overline{\eta_{1}}+ \eta_{2} \wedge \overline{\eta_{2}}) \wedge \rho \wedge R_k\right).
        \end{aligned}
    \end{equation}

    On the one hand, since $\psi_k\geq 0$ and $\psi_k=0$ on $\cup_{j=1}^m\{u_j\leq -k\}$, we have
    \begin{equation}
    \label{closdnesseq1}
        \begin{aligned}
            \int_{U} \tau^2d\psi_k \wedge d^c \psi_k \wedge \rho \wedge R_k  &\leq \frac{1}{2} \int_{U} \tau^2 \ddc\psi_{k}^{2} \wedge \rho \wedge R_{k}\\
            &=\frac{1}{2} \int_U \psi_{k}^{2} \ddc \tau^2 \wedge \rho \wedge R_{k}\\
            &= \frac{1}{2} \int_U \psi_{k}^{2} \ddc \tau^2 \wedge \rho \wedge \mathbbm{1}_{\cap_{j=1}^{m}\{u_{j}>-k\}} R_{k} \\
            &\leq c_{\rho,\tau} \lvert \mathbbm{1}_{\cap_{j=1}^{m}\{u_{j}>-k\}} R_{k} \rvert_{\supp\tau}\\
            &\leq c_{\rho,\tau} c_{n,p} \| \mathbbm{1}_{\cap_{j=1}^{m}\{u_{j}>-k\}} R_{k} \|_{\supp\tau},
        \end{aligned}
    \end{equation}
    where $c_{\rho,\tau}>0$ is a constant that depends only on $\rho$ and $\tau$. Hence, the right-hand side of (\ref{closdnesseq1}) is uniformly bounded in $k\in\N$.  

    On the other hand, since $g'$ is bounded and $g'(\psi_k)$ converges pointwise to $0$, we have
     \begin{equation}
    \label{closdnesseq2}
        \begin{aligned}
            &\int_U i(g'(\psi_k))^2 (\eta_{1} \wedge \overline{\eta_{1}}+ \eta_{2} \wedge \overline{\eta_{2}}) \wedge \rho \wedge R_k\\
            &=\int_U i(g'(\psi_k))^2 (\eta_{1} \wedge \overline{\eta_{1}}+ \eta_{2} \wedge \overline{\eta_{2}}) \wedge \rho \wedge \mathbbm{1}_{\cap_{j=1}^{m}\{u_{j}>-k\}}R_k\\
            &\leq \int_U i(g'(\psi_k))^2 (\eta_{1} \wedge \overline{\eta_{1}}+ \eta_{2} \wedge \overline{\eta_{2}}) \wedge \rho \wedge R,
        \end{aligned}
    \end{equation}
    which tends to $0$ as $k\to\infty$ by the dominated convergence theorem.

    Combining (\ref{closdnesseq0}), (\ref{closdnesseq1}), and (\ref{closdnesseq2}) completes the proof.
\end{proof}

\section{Relative non-pluripolar product on compact Hermitian manifolds}
\label{Section:CompactHermitianManifolds}

Let $(X,\omega)$ be an $n$-dimensional compact Hermitian manifold such that $\omega$ satisfies (\ref{maincondition}), namely
\[
    \partial \overline{\partial} \omega = 0,\quad
    \partial \omega \wedge \overline{\partial} \omega = 0.
\]
In this section, we show that the relative non-pluripolar product is well-defined on $X$ and satisfies the monotonicity property. This generalizes Vu's result in \cite{Vu_RelativeNonPluripolar} from the K\"ahler setting to the Hermitian setting. We first establish some equivalent conditions for (\ref{maincondition}).

\begin{lemma}
\label{omegacondition}
    Let $(X,\omega)$ be a Hermitian manifold. The following statements are equivalent.
    \begin{itemize}
        \item[\upshape(i)]
        $\partial \overline{\partial} \omega=0$ and $\partial \omega \wedge \overline{\partial}\omega=0$.

        \item[\upshape(ii)]
        $\partial \overline{\partial} \omega^k=0$ for every $k\in\N$.
    \end{itemize}
\end{lemma}

\begin{proof}
    For every $k\in\N$, we have
    \[
        \partial \overline{\partial} \omega^{k+1}=\omega\wedge\partial \overline{\partial} \omega^k+\omega^k\wedge\partial \overline{\partial} \omega+2k\omega^{k-1}\wedge\partial\omega\wedge\overline{\partial}\omega.
    \]
    Therefore,  (i) and (ii) are equivalent.
\end{proof}

\begin{theorem}
\label{Theorem:WellDefined}
    Let $(X,\omega)$ be an $n$-dimensional compact Hermitian manifold such that $\omega$ satisfies (\ref{maincondition}). Let $T_j$ be a closed positive $(1,1)$-current on $X$ for $j=1,\dots,m$. Let $T$ be a closed positive $(p,p)$-current on $X$ such that $m+p \leq n$. Then the relative non-pluripolar product $\langle T_{1} \wedge \dots \wedge T_{m} ~\dot{\wedge}~ T \rangle$ is well-defined.
\end{theorem}

\begin{proof}
    For every $j=1,\dots,m$, there exists a real closed $(1,1)$-form $\theta_{j}$ on $X$ and a quasi-psh function $\tilde{u}_j$ on $X$ such that
    \[
        T_j=dd^c \tilde{u}_j +\theta_{j}.
    \]
    Since $X$ is compact, there exists a constant $c>0$ such that $\theta_j\leq c\omega$ for every $j=1,\dots,m$. For every $k\in\N$, set
    \[
	    \tilde{u}_{j,k}:=\max\{\tilde{u}_j, -k\},\quad \tilde{R}_k:=(dd^c \tilde{u}_{1,k}+\theta_{1})\wedge\dots\wedge (dd^c \tilde{u}_{m,k}+\theta_{m})\wedge T.
	\]
	It follows from Lemma \ref{Lemma:UniformlyBounded_QuasiPsh} (i) that
    \[
        0\leq\mathbbm{1}_{\cap_{j=1}^m \{\tilde{u}_j>-k\}}\tilde{R}_k\leq\mathbbm{1}_{\cap_{j=1}^m \{\tilde{u}_j>-k\}}(dd^c\tilde{u}_{1,k}+c\omega)\wedge\dots\wedge(dd^c \tilde{u}_{m,k}+c\omega)\wedge T.
    \]
    
    We now prove that $(dd^c \tilde{u}_{1,k}+c\omega)\wedge\dots\wedge(dd^c \tilde{u}_{m,k}+c\omega)\wedge T$ is a positive current. Since
    \begin{align*}
        dd^c \tilde{u}_j+c\omega\geq dd^c \tilde{u}_j+\theta_j\geq0, \quad dd^c (-k) +c\omega= c\omega\geq 0,
    \end{align*}
    $\tilde{u}_{j,k}\in PSH(X,c\omega)$ for every $j$ and $k$. By \cite[Theorem 1]{Blocki-Kolodziej_RugularizationPshManifolds}, there exists a sequence
    \[
        (\tilde{u}^\ell_{j,k})_{\ell\in\N}\subset\mathrm{PSH}(X, c\omega)\cap\mathscr{C}^{\infty}(X,\R)
    \]
    decreases to $\tilde{u}_{j,k}$ as $\ell\to\infty$. It follows that the current
    \[
        (dd^c \tilde{u}^\ell_{1,k}+c\omega)\wedge\dots\wedge(dd^c \tilde{u}^\ell_{m,k}+c\omega)\wedge T
    \]
    is positive for every $\ell\in\N$ and weakly converges to
    \[
        (dd^c \tilde{u}_{1,k}+c\omega)\wedge\dots\wedge(dd^c \tilde{u}_{m,k}+c\omega)\wedge T
    \]
    as $\ell\to\infty$. Hence, $(dd^c \tilde{u}_{1,k}+c\omega)\wedge\dots\wedge(dd^c \tilde{u}_{m,k}+c\omega)\wedge T$ is a positive current on $X$.

    It follows from Lemma \ref{omegacondition} that $dd^c \omega^k=0$ for every $k\in\N$, and hence we have
    \begin{align*}
        \| \mathbbm{1}_{\cap_{j=1}^m \{\tilde{u}_j>-k\}}\tilde{R}_k\|&\leq\lVert\mathbbm{1}_{\cap_{j=1}^m \{\tilde{u}_j>-k\}}(dd^c \tilde{u}_{1,k}+c\omega)\wedge\dots\wedge(dd^c \tilde{u}_{m,k}+c\omega)\wedge T\rVert\\
        &\leq \|(dd^c \tilde{u}_{1,k}+c\omega)\wedge\dots\wedge(dd^c \tilde{u}_{m,k}+c\omega)\wedge T\|\\
        &=\int_{X}(dd^c \tilde{u}_{1,k}+c\omega)\wedge\dots\wedge(dd^c \tilde{u}_{m,k}+c\omega)\wedge T\wedge\omega^{n-m-p}\\
        &=\lim_{\ell\to\infty}\int_X (dd^c \tilde{u}^\ell_{1,k}+c\omega)\wedge\dots\wedge(dd^c \tilde{u}^\ell_{m,k}+c\omega)\wedge T\wedge\omega^{n-m-p}\\
        &=c^m \lVert T\rVert.
    \end{align*}
    Therefore, the relative non-pluripolar product $\langle T_{1} \wedge \dots \wedge T_{m} \ \dot{\wedge} \ T \rangle$ is well-defined. 
\end{proof}

All compact Riemann surfaces are K\"ahler. Every $n$-dimensional compact Hermitian manifold admits a Gauduchon form, i.e., a Hermitian form $\omega$ satisfying $dd^c \omega^{n-1}=0$. In particular, on every compact complex surface, a Gauduchon form clearly satisfies condition (ii) of Lemma \ref{omegacondition} and hence (\ref{maincondition}). Inspired by \cite{Tosatti-Weinkove_EstimatesComplexMA}, we present the following examples of compact non-K\"ahler manifolds that admit Hermitian forms satisfying (\ref{maincondition}).
     
\begin{example}
\label{Example:nonKahler}
    Let $Y$ be a Hopf surface, which is a compact non-K\"ahler manifold, and let $\omega_Y$ be a Gauduchon form on $Y$. Let $Z$ be an $n$-dimensional compact K\"{a}hler manifold with K\"{a}hler form $\omega_Z$. Then
    \[
        X:=Y\times Z
    \]
    is an $(n+2)$-dimensional compact non-K\"{a}hler manifold. Denote by
    \[
        p_Y: X\to Y,\quad p_Z:X\to Z
    \]
    the natural projections. Then the Hermitian form
    \[
        \omega_X:=p^*_Y\ \omega_Y +p^*_Z\ \omega_Z
    \]
    on $X$ satisfies (\ref{maincondition}). 
\end{example}

Let $(X, \omega ) $ be an $n$-dimensional compact Hermitian manifold. The \textit{upper volume} of $\omega$ is defined by
\[
    \overline{\vol} (X,\omega) : = \sup \left\{ \int_{X} \left( \ddc u + \omega\right)^{n} : u \in \PSH(X,\omega) \cap \mathscr{C}^{\infty}(X,\R) \right\}\in [0,\infty],
\]
When $\overline{\vol}(X,\omega) < \infty$, we say that $X$ satisfies the \textit{bounded volume condition}. We note that the bounded volume condition is independent of the choice of the Hermitian form $\omega$, and until now, it is still unclear whether every compact Hermitian manifold of dimension at least $3$ satisfies the bounded volume condition. The non-pluripolar product is always well-defined if $X$ satisfies the bounded volume condition (\cite[Theorem A]{BottChernvolumes}).

Let $V$ be an irreducible analytic subset of $X$ of dimension $n-p$. We denote by $[V]$ the current of integration along the regular part of $V$. It is still unclear whether the bounded volume condition guaranties the well-definedness of the non-pluripolar product relative to an arbitrary closed positive $(p,p)$-current $T$. The following remark shows, however, that under the bounded volume condition, the non-pluripolar product relative to $[V]$ is always well-defined.

\begin{remark}
\label{bvcrelativenon-pluripolar}
    Let $(X,\omega)$ be an $n$-dimensional compact Hermitian manifold such that $X$ satisfies the bounded volume condition. Let $V$ be an irreducible analytic subset of $X$ of dimension $n-p$. Let $T_{1} , \dots , T_{m}$ be closed positive $(1,1)$-currents such that $m+ p \leq n$. Then the relative non-pluripolar product $\langle T_{1} \wedge \dotsi \wedge T_{m} ~\dot{\wedge}~[V] \rangle$ is well-defined.
\end{remark}
\begin{proof}
    By Hironaka's desingularization theorem, there exists a modification $\sigma \colon X' \rightarrow X$, given by a composition of finitely many blow-ups along smooth centers, such that the strict transform $V'$ is a closed submanifold of $X'$ of dimension $n-p$. Since $\sigma$ is bimeromorphic, it follows from \cite[Theorem A]{Quasiplurienvelopes2} that $X'$ satisfies the bounded volume condition. Set $\omega' := \sigma^{*}\omega$, which is a semi-positive form on $X'$. Let $\omega_{X'}$ be a Hermitian form on $X'$ such that $\omega_{X'}\geq \omega'$, and $\iota : V' \hookrightarrow X'$ be the natural inclusion. By \cite[Theorem A]{Plurisignedhermitianmetrics}, we have $\overline{\vol}(V',\iota^*\omega_{X'})<  \infty$.

    For every $u \in \PSH(X,\omega) \cap \mathscr{C}^{\infty}(X,\R)$, set $u': = \sigma^{*} u$. Then $ u'\in \PSH(X',\omega')$ and 
    \begin{align*}
        \int_{X} (\ddc u + \omega)^{n-p}\wedge[V] = \int_{X'}  (\ddc u' + \omega' )^{n-p}\wedge[V'].
    \end{align*}
    Since $\omega_{X'}\geq \omega'$, we have $u'\in PSH (X',\omega_{X'})$, $\iota^*u'\in PSH (V',\iota^*\omega_{V'})\cap \mathscr{C}^\infty(V',\R)$, and
    \begin{align}\label{remark3.4ineq0}
        \begin{aligned}
            \int_{X} (\ddc u + \omega)^{n-p}\wedge [V]&= \int_{X'} (\ddc u' + \omega' )^{n-p}\wedge [V']\\
        &\leq \int_{X'}(\ddc u' + \omega_{X'} )^{n-p}\wedge [V']\\
        &=\int_{V'} (\ddc \iota^*u' + \iota^*\omega_{X'} )^{n-p}\\
        &\leq \overline{\vol}(V', \iota^*\omega_{X'})\\
        &<  \infty.
        \end{aligned}
    \end{align}
    For every $u\in PSH (X,\omega)\cap L^\infty(X)$, after taking a decreasing approximation sequence in $PSH (X,\omega)\cap \mathscr{C}^\infty(X,\R)$ and applying (\ref{remark3.4ineq0}), we deduce that
    \begin{equation}
    \label{remark3.4ineq}
        \int_{X} (\ddc u + \omega)^{n-p}\wedge[V] \leq \overline{\vol}(V', \iota^*\omega_{X'})<  \infty.
    \end{equation}
    
    We retain the notation from the proof of Theorem \ref{Theorem:WellDefined}, except that $T$ is replaced by $[V]$. It follows from (\ref{remark3.4ineq}) that
    \begin{align*}
        \| \mathbbm{1}_{\cap_{j=1}^m \{\tilde{u}_j>-k\}}\tilde{R}_k\|&\leq\lVert\mathbbm{1}_{\cap_{j=1}^m \{\tilde{u}_j>-k\}}(dd^c \tilde{u}_{1,k}+c\omega)\wedge\dots\wedge(dd^c \tilde{u}_{m,k}+c\omega)\wedge [V]\rVert\\
        &\leq \int_{X}(dd^c \tilde{u}_{1,k}+c\omega)\wedge\dots\wedge(dd^c \tilde{u}_{m,k}+c\omega)\wedge [V]\wedge\omega^{n-m-p}\\
        &\leq \int_X \big( dd^c(\tilde{u}_{1,k}+\dots +\tilde{u}_{m,k})+(cm+n-m-p)\omega \big)^{n-p}\wedge[V] \\
        &\leq \frac{\overline{\vol}(V',\iota^*\omega_{X'})}{(cm+n-m-p)^{n-p}},
    \end{align*}
    which is independent of $k\in\N$ and finite, therefore $\langle T_1\wedge\cdots\wedge T_m~\dot{\wedge}~[V]\rangle$ is well-defined.
\end{proof}

We now aim to prove the monotonicity property of relative non-pluripolar products on compact Hermitian manifolds, stated in the following.

\begin{theorem}
\label{hermitianmonotonicity}
    Let $(X,\omega)$ be an $n$-dimensional compact Hermitian manifold such that $\omega$ satisfies (\ref{maincondition}). For $j=1,\dots, m$, let $\{\theta_{j}\}_{BC}$ be a pseudoeffective Bott--Chern class, and $T_{j},T'_{j}  \in \{\theta_{j}\}_{BC}$ be closed positive $(1,1)$-currents such that $T'_{j}$ is less singular than $T_{j}$. Let $T$ be a closed positive $(p,p)$-current on $X$ such that $m+ p \leq n$. Then we have 
    \[
        \int_{X}\langle T_{1} \wedge \dots \wedge T_{m} ~\dot{\wedge}~ T \rangle \wedge \omega^{n-m-p}  \leq \int_{X}\langle T'_{1} \wedge \dots \wedge T'_{m} ~\dot{\wedge}~ T \rangle \wedge \omega^{n-m-p}.
    \]
\end{theorem}
    
We review some results from \cite{Vu_RelativeNonPluripolar} that will be used in the proof of Theorem \ref{hermitianmonotonicity}. For the reader's convenience, we include the proofs with more details.

\begin{lemma} (\cite[Lemma 4.1]{Vu_RelativeNonPluripolar})
\label{localmonoton}
    Let $U \subset \C^{n}$ be an open subset. For $j=1,\dots,m$, let $u_{j}\in PSH(U)$ and $(u_{j}^{\ell})_{\ell\in\N}\subset PSH (U)$, such that $u^{\ell}_{j} \geq u_{j}$ and $u^{\ell}_{j} \to u_{j}$ in $L^{1}_{loc}$ as $\ell \to \infty$. Let $T$ be a closed positive $(p,p)$-current on $U$ such that $m+p \leq n$. Assume that the non-pluripolar products
    \[
        \langle \ddc u_1\wedge\dots\wedge dd^c u_m ~\dot{\wedge}~ T \rangle,\quad \langle \ddc u_1^\ell\wedge\dots\wedge dd^c u_m^\ell ~\dot{\wedge}~ T \rangle
    \]
    are well-defined for every $\ell\in\N$. Then, we have
    \[
        \int_{U} \langle \ddc u_1\wedge\dots\wedge dd^c u_m ~\dot{\wedge}~ T \rangle \wedge \Phi\leq\liminf_{\ell \to \infty}\int_{U} \langle \ddc u_1^\ell\wedge\dots\wedge dd^c u_m^\ell ~\dot{\wedge}~ T \rangle \wedge \Phi
    \]
    for every strongly positive smooth $(n-m-p,n-m-p)$-form $\Phi$ on $U$ with compact support.
\end{lemma}

\begin{proof}
    Let $\Phi$ be a strongly positive smooth $(n-m-p,n-m-p)$-form on $U$ with compact support. By Hartogs' lemma, $u_j^\ell$ is locally uniformly bounded from above in $j$ and $\ell$. Hence, by shrinking $U$, we can assume that there exists a constant $c\geq 1$ such that
    \begin{align*}
        |z|\leq c, \quad u_j\leq u_j^\ell \leq  c
    \end{align*}
    on $U$ for all $j$ and $\ell$. For every $k\in \N$, set
    \begin{align*}
        &u_{j,k} := \max\{u_{j},-k\},\quad&R_{k}:= \wedge_{j=1}^{m} \ddc u_{j,k} \wedge T,\qquad&R:= \langle \wedge_{j=1}^{m} \ddc u_{j} ~\dot{\wedge}~T \rangle,\\
        &u^{\ell}_{j,k} := \max \{u^{\ell}_{j},-k\},\quad&R^{\ell}_{k}:= \wedge_{j=1}^{m} \ddc u^{\ell}_{j,k} \wedge T,\qquad&R^{\ell}:= \langle \wedge_{j=1}^{m} \ddc u_{j}^{\ell} ~\dot{\wedge}~T \rangle.
    \end{align*}
    It follows from $u^{\ell}_{j} \geq u_{j}$ and Lemma \ref{Lemma:UniformlyBounded_Psh}, that $\{u_{j,k}>-k\}\subset\{u_j^\ell>-k\}$ and
    \begin{align}
    \label{localmonotonequ1}
        \int_{U} R^{\ell} \wedge \Phi \geq  \int_{U} \mathbbm{1}_{\cap_{j=1}^{m} \{u_{j}^{\ell}>-k\}} R^{\ell}_{k} \wedge \Phi\geq \int_{U} \mathbbm{1}_{\cap_{j=1}^{m} \{u_{j,k}>-k\}}  R^{\ell}_{k}  \wedge \Phi.
    \end{align}
    
    Fix an arbitrary $k\in\N$ such that $k\geq c^2$. Since $u_{j,k}$ is bounded on $U$ for $j=1,\dots,m$, by \cite[Theorem 2.4]{Vu_RelativeNonPluripolar}, for every $\varepsilon>0$ there exists an open subset $U_{\varepsilon}$ of $U$ such that $\capacity_{T}(U_{\varepsilon},U) < \varepsilon$ and $u_{j,k}|_{U \backslash U_{\varepsilon}}$ is continuous for every $j$. By Tietze extension theorem, we can extend $u_{j,k}|_{U \backslash U_{\varepsilon}}$ to a function $\tilde{u}_{j,k}\in \mathscr{C}(U,\R)$ for every $j$.

    Since $u_{j,k}, (u^{\ell}_{j,k})_{\ell\in\N}$ are uniformly bounded on $U$, $u^{\ell}_{j,k} \geq u_{j,k}$ and $u_{j,k}^{\ell} \to u_{j,k}$ in $L_{loc}^{1}$ as $\ell \to \infty$, by \cite[Theorem 2.2]{Vu_RelativeNonPluripolar}, we have $R^{\ell}_{k} \to R_{k}$  as $\ell \to \infty$. It then follows that
    \begin{equation}
    \label{localmonotonequ2}
        \int_{U}\mathbbm{1}_{\cap_{j=1}^{m}\{\tilde{u}_{j,k}>-k\}} R_{k} \wedge \Phi\leq\liminf_{\ell \to\infty} \int_{U} \mathbbm{1}_{\cap_{j=1}^{m}\{\tilde{u}_{j,k}>-k\}}R^{\ell}_{k} \wedge \Phi, 
    \end{equation}
    since $\mathbbm{1}_{\cap_{j=1}^{m}\{\tilde{u}_{j,k}>-k\}}$ is lower semi-continuous.
    
    Notice that $|z|^2\leq c^2\leq k$ on $U$ and $-k\leq u_{j,k}\leq u_{j,k}^{\ell}\leq c\leq c^2 \leq k$ on $U$, then we have
    \begin{equation}
    \label{localmonotonequ3}
        \begin{aligned}
            &\Big{|}\int_{U} \mathbbm{1}_{\cap_{j=1}^{m}\{\tilde{u}_{j,k}>-k\}}R^{\ell}_{k} \wedge \Phi  -\int_{U} \mathbbm{1}_{\cap_{j=1}^{m}\{u_{j,k}>-k\}}R^{\ell}_{k} \wedge \Phi \Big{|}\\
            &\leq 2 c_\Phi \lvert \mathbbm{1}_{U_\varepsilon}R^{\ell}_{k}\rvert_{\supp\Phi}\\
            &\leq 2c_\Phi c_{n,p} \int_{\supp\Phi}\mathbbm{1}_{U_{\varepsilon}} R_{k}^\ell \wedge (dd^c|z|^2)^{n-m-p}\\
            &\leq 2 (2k)^{n-p}c_\Phi c_{n,p} \capacity_{T}(U_{\varepsilon},U)\\
            &< 2 (2k)^{n-p}c_\Phi c_{n,p}  \varepsilon,
        \end{aligned}
    \end{equation}
    where $c_\Phi>0$ is a constant that depends only on $\Phi$. Similarly, we have
    \begin{equation}
    \label{localmonotonequ4}
        \begin{aligned}
            &\Big{|}\int_{U} \mathbbm{1}_{\cap_{j=1}^{m}\{\tilde{u}_{j,k}>-k\}}R_{k} \wedge \Phi  -\int_{U} \mathbbm{1}_{\cap_{j=1}^{m}\{u_{j,k}>-k\}}R_{k} \wedge \Phi \Big{|}< 2 (2k)^{n-p}c_\Phi c_{n,p}  \varepsilon.
        \end{aligned}
    \end{equation}
    
    It follows from (\ref{localmonotonequ1}), (\ref{localmonotonequ2}), (\ref{localmonotonequ3}), and (\ref{localmonotonequ4}) that 
    \[
        \liminf_{\ell \to \infty} \int_{U}R^{\ell} \wedge \Phi \geq \int_{U} \mathbbm{1}_{\cap_{j=1}^{m}\{u_{j,k}>-k\}}R_{k} \wedge \Phi  -4(2k)^{n-p}c_\Phi c_{n,p}\epsilon.
    \]
    By letting $\varepsilon\to 0$ and then $k\to\infty$, we obtain the desired inequality.
\end{proof}

We now turn to the global setting. Let $T$ and $T'$ be closed positive $(1,1)$-currents in the same Bott--Chern class. We say that $T$ and $T'$ are of the same singularity type if each is less singular than the other. We then have the following proposition.

\begin{proposition} (\cite[Proposition 4.2]{Vu_RelativeNonPluripolar})
\label{samesingularities}
    Let $(X,\omega)$ be an $n$-dimensional compact Hermitian manifold. For $j=1,\dots,m$, let $T_{j},T_{j}'$ be two closed positive $(1,1)$-currents  in the same Bott--Chern class and of the same singularity type. Let $T$ be a closed positive $(p,p)$-current such that $m + p \leq n$. Assume that for every $J , J' \subset \{1, \dots , m\}$ such that $J \cap J' = \emptyset$, the non-pluripolar product 
    \[
        \left\langle \bigwedge_{j \in J} T_{j} \wedge \bigwedge_{j' \in J'} T'_{j'} ~\dot{\wedge}~ T \right\rangle
    \]
    is well-defined. Then, for every $\ddc$-closed $(n-m-p,n-m-p)$-form $\Phi$ on $X$, we have 
    \[
        \int_{X} \langle T_1\wedge\dots\wedge T_m ~\dot{\wedge}~ T \rangle \wedge \Phi = \int_{X} \langle T'_1\wedge\dots\wedge T'_m ~\dot{\wedge}~ T \rangle \wedge \Phi.
    \]
\end{proposition}

\begin{proof}
    For every $j=1,\dots,m$, write
    \[
        T_{j} = \ddc \tilde{u}_j + \theta_{j},\quad T'_{j} = \ddc \tilde{u}'_j + \theta_{j}
    \]
    where $\theta_j$ is a real closed $(1,1)$-form on $X$ and $\tilde{u}_j,\tilde{u}'_j\in PSH (X,\theta_j)$ are negative functions. For every $k\in\N$ and $J , J' \subset \{1, \dots , m\}$ such that $J \cap J' = \emptyset$, set
    \begin{align*}
        &\tilde{u}_{j,k}:=\max\{\tilde{u}_j,-k\},\quad \tilde{u}'_{j,k}:=\max\{\tilde{u}'_j,-k\},\quad w_{j,k}:= \tilde{u}_{j,k} - \tilde{u}'_{j,k},\\
        &\psi_{k} := \frac{1}{k} \max \left\{\sum_{j=1}^{m}(\tilde{u}_{j}+\tilde{u}'_{j}),-k \right\}+1,\quad A:= \bigcup_{j=1}^m\{\tilde{u}_j=-\infty\}=\bigcup_{j=1}^m\{\tilde{u}'_j=-\infty\},\\
         &R_{J,J',k}:=  \bigwedge_{j \in J}(dd^c \tilde{u}_{j,k}+\theta_j) \wedge \bigwedge_{j' \in J'} (dd^c \tilde{u}'_{j,k}+\theta_j) \wedge T,\quad R_{J,J'}:=\left\langle \bigwedge_{j \in J} T_{j} \wedge \bigwedge_{j' \in J'} T'_{j'} ~\dot{\wedge}~ T \right\rangle.
    \end{align*}
    Since $T_{j},T'_{j}$ are of the same singularity type, $w_{j,k}$ is uniformly bounded in $j$ and $k$. Notice that $\psi_{k}$ is a quasi-psh function on $X$, $0\leq \psi_k\leq 1$, $\psi_k=0$ on $\cup_{j=1}^m\{\tilde{u}_j\leq -k\}\cup\{\tilde{u}'_j\leq -k\}$, and $\psi_{k}$ increases to $\mathbbm{1}_{X \backslash A}$ as $k\to\infty$. Set
    \[
    S := \langle T_{1}\wedge \cdots\wedge T_m ~\dot{\wedge}~ T \rangle - \langle T'_{1}\wedge \cdots\wedge T'_m ~\dot{\wedge}~T \rangle.
    \]
    Then $\psi_kS$ weakly converges to $\mathbbm{1}_{X\backslash A}S$, which is equal to $S$ by Lemma \ref{Lemma:UniformlyBounded_QuasiPsh} (iii). Hence, the desired assertion is equivalent to
    \begin{equation}
    \label{samesingularitieseq1}
        \lim_{k\to\infty}\int_{X}\psi_{k} S \wedge \Phi =0
    \end{equation}
    for every $\ddc$-closed $(n-m-p,n-m-p)$-form $\Phi$  on $X$.
    By the construction of $\psi_{k}$ and Lemma \ref{Lemma:UniformlyBounded_QuasiPsh} (iii), we have
    \begin{align*}
        \psi_{k} S&=\psi_k \mathbbm{1}_{\cap_{j=1}^m \{\tilde{u}_j>-k\}}\langle T_{1}\wedge \cdots\wedge T_m ~\dot{\wedge}~ T \rangle - \psi_k \mathbbm{1}_{\cap_{j=1}^m \{\tilde{u}'_j>-k\}}\langle T'_{1}\wedge \cdots\wedge T'_m ~\dot{\wedge}~ T \rangle\\
        &= \psi_{k} \big((dd^c \tilde{u}_{1,k}+\theta_1)\wedge\dots\wedge(dd^c \tilde{u}_{m,k}+\theta_1) -  (dd^c \tilde{u}'_{1,k}+\theta_1)\wedge\dots\wedge(dd^c \tilde{u}'_{m,k}+\theta_1)\big)\wedge T\\
        &=\sum_{s=1}^m \psi_k dd^c w_{s,k}\wedge R_{\{1,\dots,s-1\},\{s+1,\dots,m\},k}.
    \end{align*}
    Therefore, it suffices to show that, for every $s\notin J\cup J'$ and $\ddc$-closed $(r,r)$-form $\Phi$  on $X$,
    \[
        \lim_{k\to\infty}\int_{X} \psi_{k} \ddc w_{s,k} \wedge R_{J,J',k} \wedge \Phi=0,
    \]
    here $r:=n-|J|-|J'|-p-1$.

    Let $A'>0$ be a constant such that $A'\omega\geq \theta_j$ for $j=1,\dots,m$, then for every $j$ and $k$,
    \[
        \tilde{u}_{j,k},\tilde{u}'_{j,k}\in PSH (X,A'\omega),\quad \psi_k\in PSH (X, \frac{2mA'}{k}\omega)\subset PSH (X,2mA'\omega).
    \]
    By \cite[Theorem 1]{Blocki-Kolodziej_RugularizationPshManifolds}, there exists sequences of smooth functions in $PSH (X,A'\omega)$ decreasing to $\tilde{u}_{j,k}$ and $\tilde{u}'_{j,k}$, respectively, and a sequence of smooth functions in $PSH (X,\frac{2mA'}{k}\omega)$ decreasing to $\psi_k$.

    Although $R_{J,J',k}$ is locally the difference of two closed positive currents, it is not necessarily positive. Hence we cannot estimate its mass directly and we modify $R_{J,J',k}$ as follows. Set $\tilde{R}_{J,J'} := \mathbbm{1}_{X \backslash A} R_{J,J'}$. Then $\tilde{R}_{J,J'}$ is positive and closed since $A$ is a complete pluripolar subset of $X$ (see \cite[Remark 1.10]{Boucksom-Eyssidieux-Guedj-Zeriahi}). Moreover, by the construction of $\psi_{k}$ and Lemma \ref{Lemma:UniformlyBounded_QuasiPsh} (iii), we have
    \begin{align*}
        \psi_{k} \tilde{R}_{J,J'} &= \psi_{k} R_{J,J'}\\
        &=\psi_k\mathbbm{1}_{\cap_{j\in J}\{\tilde{u}_j>-k\}}\mathbbm{1}_{\cap_{j'\in J'}\{\tilde{u}'_j>-k\}}R_{J,J'} \\
        &=\psi_k\mathbbm{1}_{\cap_{j\in J}\{\tilde{u}_j>-k\}}\mathbbm{1}_{\cap_{j'\in J'}\{\tilde{u}'_j>-k\}}R_{J,J',k} \\
        &=\psi_{k} R_{J,J',k}.
    \end{align*}
    Next, applying integration by parts to smooth decreasing approximations of $\tilde{u}_{s,k}$ and $\tilde{u}'_{s,k}$, then passing to the limit, and using the fact that $dd^c\Phi=0$, we have
    \begin{align*}
        \int_{X} \psi_{k} \ddc w_{s,k} \wedge R_{J,J',k} \wedge \Phi =  2\int_{X} w_{s,k} d \psi_{k} \wedge \dc \Phi \wedge\tilde{R}_{J,J'} + \int_{X} w_{s,k}  \ddc \psi_{k} \wedge \tilde{R}_{J,J'} \wedge \Phi.
    \end{align*}
    
    We estimate the first term on the right-hand side. $\dc \Phi$ can be written as a finite sum of terms of the form $\tau\wedge\Psi$, where every $\tau$ is a smooth form on $X$ and is of bi-degree $(1,0)$ or $(0,1)$, every $\Psi$ is a strongly positive smooth $(r,r)$-form on $X$ and $\Psi\leq c_\Phi \omega^r$ for some constant $c_\Phi>0$. Hence, we can assume that $\dc \Phi = \tau_{1} \wedge \Psi_{1}+ \overline{\tau_{2}} \wedge \Psi_{2}$, where $\tau_{1}, \tau_{2}$ are $(1,0)$-forms and $\Psi_{1},\Psi_{2}$ are strongly positive. Applying the Cauchy–Schwarz inequality to the smooth decreasing approximations of $\psi_k, \tilde{u}_{s,k},\tilde{u}'_{s,k}$, then passing to the limit, we have
    \begin{align*}
        &\Big{|}\int_{X}  w_{s,k} d \psi_{k} \wedge \dc \Phi \wedge \tilde{R}_{J,J'} \Big{|}^{2}\\ 
        &\lesssim  \Big{(} \int_{X}  w_{s,k}^{2} d \psi_{k} \wedge d^c \psi_{k} \wedge \tilde{R}_{J,J'} \wedge \Psi_{1} \Big{)}\Big{(} \int_{X} i \tau_{1} \wedge \overline{\tau_{1}}  \wedge \tilde{R}_{J,J'} \wedge \Psi_{1}  \Big{)}\\
        &+ \Big{(} \int_{X}  w_{s,k}^{2} d \psi_{k} \wedge d^c \psi_{k} \wedge \tilde{R}_{J,J'} \wedge \Psi_{2} \Big{)}\Big{(} \int_{X} i \tau_{2} \wedge \overline{\tau_{2}} \wedge \tilde{R}_{J,J'} \wedge \Psi_{2}  \Big{)}\\
        &\lesssim \Big{(} \int_{X}  d \psi_{k} \wedge \dc \psi_{k} \wedge \tilde{R}_{J,J'} \wedge \omega^r\Big{)} \|R_{J,J'}\|.
        \end{align*}
    Then, it suffices to show that $d\psi_k\wedge d^c\psi_k\wedge \tilde{R}_{J,J'}$ weakly converges to $0$ as $k\to\infty$. On the one hand, we have $\psi_k\in PSH (X, 2mA'\omega)$ for every $k\in \N$. For every small enough local chart $U$, there exists $h\in\mathscr{C}^\infty(U,\R)$ such that $dd^c h\geq 2mA'\omega$ on $U$. Then we have $\psi_k|_U+h\in PSH (U)$, $h\in PSH(U)$, and $\psi_k|_U+h$ increases to $h$ outside the complete pluripolar set $A\cap U$. By the fact that $\tilde{R}_{J,J'}|_U$ has no mass on $A\cap U$ and \cite[Remark 2.7]{Vu_RelativeNonPluripolar},
    \begin{align*}
        &(\psi_kdd^c \psi_k\wedge\tilde{R}_{J,J'})|_U\\
        &= \Big((\psi_k|_U +h)dd^c (\psi_k|_U +h)-(\psi_k|_U +h)dd^c h- h dd^c (\psi_k|_U +h) + hdd^c h\Big)\wedge\tilde{R}_{J,J'}|_U
    \end{align*}
    and weakly converges to $0$ on $U$. It follows that $\psi_kdd^c \psi_k\wedge\tilde{R}_{J,J'}$ weakly converges to $0$ on $X$. On the other hand, by the dominated convergence theorem we have
    \begin{align*}
        \lim_{k\to\infty}dd^c\psi_k^2\wedge\tilde{R}_{J,J'}= \lim_{k\to\infty}dd^c(\psi_k^2\tilde{R}_{J,J'})= dd^c (\mathbbm{1}_{X\backslash A}\tilde{R}_{J,J'})=dd^c \tilde{R}_{J,J'}=0.
    \end{align*}
    Therefore,
    \[
        \lim_{k\to\infty}d \psi_{k} \wedge \dc \psi_{k} \wedge \tilde{R}_{J,J'} =  \lim_{k\to\infty}(\frac{1}{2}\ddc \psi_{k}^{2} - \psi_{k}\ddc \psi_{k}) \wedge \tilde{R}_{J,J'}= 0.
    \]
    
    It remains to prove that
    \[
        \lim_{k \to \infty}\int_{X} w_{s,k} \ddc \psi_{k} \wedge \tilde{R}_{J,J'}\wedge \Phi=0.
    \]
    $\Phi$ can be written as a finite linear combination of strongly positive smooth forms with constant coefficients, hence, we can assume that $\Phi$ is strongly positive. Then, we have
    \begin{align*}
        &\left| \int_{X} w_{s,k} \ddc \psi_{k} \wedge \tilde{R}_{J,J'}\wedge \Phi\right|\\
        &=\left| \int_{X} w_{s,k} (\ddc \psi_{k}+\frac{2mA'}{k}\omega) \wedge \tilde{R}_{J,J'}\wedge \Phi-\frac{2mA'}{k}\int_{X} w_{s,k} ~\omega \wedge \tilde{R}_{J,J'}\wedge \Phi\right|\\
        &\lesssim \int_X (\ddc \psi_{k}+\frac{2mA'}{k}\omega) \wedge \tilde{R}_{J,J'}\wedge \Phi + \frac{2mA'}{k}\int_{X} \omega \wedge \tilde{R}_{J,J'}\wedge \Phi,
    \end{align*}
    which tends to $0$ as $k\to\infty$, since 
    \[
        \lim_{k\to\infty}dd^c \psi_k\wedge \tilde{R}_{J,J'}= \lim_{k\to\infty}dd^c(\psi_k\tilde{R}_{J,J'})= dd^c (\mathbbm{1}_{X\backslash A}\tilde{R}_{J,J'})=dd^c \tilde{R}_{J,J'}=0.
    \]
    
    This finishes the proof.
\end{proof}

We now prove Theorem \ref{hermitianmonotonicity}, following the argument in \cite[Theorem 4.4]{Vu_RelativeNonPluripolar}.

\begin{proof}[Proof of Theorem \ref{hermitianmonotonicity}]
    By Theorem \ref{Theorem:WellDefined}, the relative non-pluripolar products are always well-defined on $X$. For every $j=1,\dots,m$, write 
    \[
        T_{j} = \ddc \tilde{u}_{j} + \theta_{j},\quad T'_{j} =\ddc \tilde{u}'_{j} + \theta_{j},
    \]
    where $\tilde{u}_{j},\tilde{u}'_{j} \in \PSH(X,\theta_{j})$. Foe every $\ell\in\N$, set
    \[
        \tilde{u}_j^\ell:= \max \{\tilde{u}_{j},\tilde{u}'_{j}-\ell\},\quad T_j^\ell := \ddc \tilde{u}_j^\ell+\theta_{j}.
    \]
    Then $\tilde{u}_j^\ell\in PSH (X,\theta_{j})$ and $\tilde{u}_j^\ell$ decreases to $\tilde{u}_{j}$ as $\ell \to \infty$. Since $T'_{j}$ is less singular than $T_{j}$, we have that $T_j^\ell$ and $T'_{j}$ are of the same singularity type. By Proposition \ref{samesingularities}, for every $dd^c$-closed $(n-m-p,n-m-p)$-form $\Phi$ on $X$, we have
    \begin{equation}
    \label{equationb}
        \int_{X}\langle T_1^\ell \wedge \dots \wedge T_m^\ell ~\dot{\wedge}~ T \rangle \wedge \Phi = \int_{X} \langle T'_{1} \wedge \dots \wedge T'_{m} ~\dot{\wedge}~ T \rangle  \wedge \Phi, \quad \forall \ell\in\N.
    \end{equation}
    In particular, by Lemma \ref{omegacondition}, we have
    \begin{equation}
    \label{hermitianmonotonequ1}
        \int_{X}\langle T_1^\ell \wedge \dots \wedge T_m^\ell ~\dot{\wedge}~ T \rangle \wedge \omega^{n-m-p} = \int_{X} \langle T'_{1} \wedge \dots \wedge T'_{m} ~\dot{\wedge}~ T \rangle  \wedge \omega^{n-m-p}, \quad \forall \ell\in\N.
    \end{equation}
    Hence, by the compactness of the space of positive currents, there exists a subsequence
    of $\big(\langle T_1^\ell \wedge \dots \wedge T_m^\ell ~\dot{\wedge}~ T \rangle\big)_{\ell\in\N}$ such that it weakly converges to a closed positive $(m+p,m+p)$-current $S$ on $X$. We denote this subsequence by $\big(\langle T_1^{r_\ell} \wedge \dots \wedge T_m^{r_\ell} ~\dot{\wedge}~ T \rangle\big)_{\ell\in\N}$.

    For every small enough local chart of $X$, there exists $h_j\in\mathscr{C}^\infty(U,\R)$ for $j=1,\dots,m$ such that $dd^c h_j=\theta_j$. Hence, we have
    \begin{align*}
        & u_j:=\tilde{u}_j|_U + h_j\in PSH(U),\quad T_j|_U=dd^c u_j,\\
        & u^l_j:=\tilde{u}_j^\ell|_U + h_j\in PSH(U),\quad T_j^\ell|_U=dd^c u^\ell_j.
    \end{align*}
    Notice that $u_j^{r_\ell}$ decrease to $u_j$ as $\ell\to\infty$, and hence $u_j^{r_\ell}\to u_j$ in $L^{1}_{loc}$. Then, it follows from Lemma \ref{localmonoton} that, for every strongly positive smooth $(n-m-p,n-m-p)$-form $\Phi$ on $U$ with compact support, we have
    \[
        \int_{U} \langle T_1\wedge\dots\wedge T_m ~\dot{\wedge}~ T \rangle \wedge \Phi\leq\liminf_{\ell \to \infty}\int_{U} \langle T_1^{r_\ell}\wedge\dots\wedge T_m^{r_\ell} ~\dot{\wedge}~ T \rangle \wedge \Phi=\int_U S\wedge\Phi.
    \]
    This implies that, for every strongly positive smooth $(n-m-p,n-m-p)$-form $\Phi$ on $X$,
    \begin{equation}
    \label{equationa}
        \int_X \langle T_1\wedge\dots\wedge T_m ~\dot{\wedge}~ T \rangle \wedge \Phi\leq\int_X S\wedge\Phi=\lim_{\ell\to\infty}\int_X \langle T_1^{r_\ell}\wedge\dots\wedge T_m^{r_\ell} ~\dot{\wedge}~ T \rangle \wedge \Phi.
    \end{equation}
    Combining (\ref{hermitianmonotonequ1}) and (\ref{equationa}), we have
    \[
        \int_{X}\langle T_{1} \wedge \dots \wedge T_{m} ~\dot{\wedge}~ T \rangle \wedge \omega^{n-m-p} \leq \int_{X} \langle T'_{1} \wedge \dots \wedge T'_{m} ~\dot{\wedge}~ T \rangle  \wedge \omega^{n-m-p}.
    \]
    This finishes the proof.
\end{proof}

Let $\alpha ,\beta \in H^{p,p}_{BC}(X,\R)$ be two $(p,p)$ Bott--Chern classes. We denote by $\alpha \geq \beta$ if there exists a closed positive $(p,p)$-current $T \in \alpha - \beta$. It follows from Lemma \ref{omegacondition} that, this relation defines a partial order on $H^{p,p}_{BC}(X,\mathbb{R})$ if $\omega$ satisfies condition \eqref{maincondition}, as in the K\"ahler case. The following remark shows that the monotonicity property (Theorem~\ref{hermitianmonotonicity}) holds in the sense of Bott--Chern cohomology classes. 

\begin{remark}
\label{monotonicityclassesremark}
    Under the assumptions in Theorem~\ref{hermitianmonotonicity}, we have
    \[
        \{\langle T_{1} \wedge \dots \wedge T_{m} ~\dot{\wedge}~ T \rangle \}_{BC} \leq \{ \langle T'_{1} \wedge \dots \wedge T'_{m} ~\dot{\wedge}~ T \rangle \}_{BC}.
    \]
\end{remark}

\begin{proof}
    We retain the notation from the proof of Theorem \ref{hermitianmonotonicity}. By (\ref{equationa}), we have
    \[
        \langle T_{1} \wedge \dots \wedge T_{m} ~\dot{\wedge}~ T \rangle\leq S,\quad \{\langle T_{1} \wedge \dots \wedge T_{m} ~\dot{\wedge}~ T \rangle \}_{BC}\leq \{S\}_{BC}.
    \]
    Hence, to obtain the desired assertion, it suffices to show that
    \[
        \{S\}_{BC}=\{\langle T'_1\wedge\dots\wedge T'_m~\dot{\wedge}~T\rangle\}_{BC}.
    \]
    By (\ref{equationb}), for every $dd^c$-closed $(n-m-p,n-m-p)$-form $\Phi$ on $X$, we have
    \begin{equation}
        \int_{X} S \wedge \Phi=\lim_{\ell\to\infty}\int_{X}\langle T_1^{r_\ell} \wedge \dots \wedge T_m^{r_\ell} ~\dot{\wedge}~ T \rangle \wedge \Phi = \int_{X} \langle T'_{1} \wedge \dots \wedge T'_{m} ~\dot{\wedge}~ T \rangle  \wedge \Phi.
    \end{equation}
    Then, by the duality between Bott--Chern cohomology groups and Aeppli cohomology groups (see \cite[Lemma 2.5]{Schweitzer} and the subsequent discussion, also \cite[Theorem 2.1]{Popovici2015AeppliGauduchon}), we have
    \[
        \{S\}_{BC}
        =
        \{\langle T'_1\wedge\cdots\wedge T'_m~\dot{\wedge}~T\rangle\bigr\}_{BC}.
    \]
    This finishes the proof.
\end{proof}

\bibliography{biblio_LiSu}
\bibliographystyle{siam}

\bigskip

\bigskip

\footnotesize
\begin{itemize}
    \item[]
    \textsc{Zhenghao Li\\Center for Complex Geometry, Institute for Basic Science,\\
    55 Expo-ro, Yuseong-gu, Daejeon 34126, Republic of Korea}\\
    \textit{E-mail address}: \texttt{zhenghaoli@ibs.re.kr}
\end{itemize}

\bigskip

\begin{itemize}
    \item[]
    \textsc{Shuang Su\\
    Center for Complex Geometry, Institute for Basic Science,\\
    55 Expo-ro, Yuseong-gu, Daejeon 34126, Republic of Korea}\\
    \textit{E-mail address}: \texttt{su30@ibs.re.kr}
\end{itemize}
\end{document}